\documentclass[12pt]{amsart}
\usepackage{amssymb}
\usepackage{amsmath}
\usepackage{mathtools}
\usepackage[all]{xy}
\usepackage{enumitem}
\usepackage{xcolor}

\renewcommand{\(}{\left(}
\renewcommand{\)}{\right)}

\newcommand{\abs}[1]{\left\lvert#1\right\rvert}
\newcommand{\norm}[1]{\left\lVert#1\right\rVert}

\newcommand{\R}{{\mathbb{R}}}

\newcommand{\e}{\varepsilon}
\renewcommand{\phi}{\varphi}
\newcommand{\p}{\rho}

\newcommand{\supp}{{\mathrm{supp}}}
\newcommand{\tr}{{\mathrm{tr}}}

\newcommand{\twisted}{\,\natural \,}

\renewcommand{\H}{{\mathcal{H}}}
\renewcommand{\S}{{\mathcal{S}}}

\theoremstyle{plain}
\newtheorem{thm}{Theorem}

\theoremstyle{definition}

\theoremstyle{remark}

\title{The Weyl transform of a compactly supported distribution}

\author{Mansi Mishra and M.~K.~Vemuri}
\address{Department of Mathematical Sciences, IIT(BHU), Varanasi 221005.}

\begin{document}

\begin{abstract}
If $T$ is a compactly supported distribution on $\R^{2n}$, then
the Weyl transform of $T$ is $p$-power traceable if and only if
the Fourier transform of $T$ is $p$-power integrable, and
the Weyl transform of $T$ is a compact operator if and only if
the Fourier transform of $T$ vanishes at infinity.
\end{abstract}


\maketitle
\thispagestyle{empty}

Let $\H = L^2\(\R^n\)$.
Let $S^p(\H)$ denote the $p$-Schatten
class of $\H$, $1\le p\le \infty$ with the convention that
$S^\infty(\H)$ is the set of bounded operators on $\H$ (see \cite{RS}).
Define 
\begin{equation*}
\(\p(x,y)\phi\)(t) = e^{\pi i (x\cdot y+2y\cdot t)}\phi(t+x),
\qquad (x,y)\in\R^{2n}, \quad \phi \in \H.
\end{equation*}
It is well known that $\p(x,y)$ is a unitary operator on $\H$
and hence $\p(x,y)\in S^\infty(\H)$ for all $(x,y)\in\R^{2n}$.

Recall that if $X \in S^1(\H)$, the {\em Fourier-Wigner transform} of $X$ is
the function on $\R^{2n}$ defined by
\begin{equation*}
\alpha(X)(x,y)= \tr(\rho(-x,-y)X)
\end{equation*}
and if $f\in L^1(\R^{2n})$, the {\em Weyl transform} of $f$ is the
operator $W(f) \in S^\infty(\H)$ defined by
\begin{equation*}
W(f) = \int_{\R^n} \int_{\R^n} f(x,y) \p(x,y) \,dx \,dy,
\end{equation*}
and that these are mutual inverses (see \cite{folland}).

The definition of the Weyl transform has been extended to tempered
distributions by several authors (see e.g., \cite{maillard,keyl});
we briefly recall this extended definition now.
If $k$ is a Schwartz class function on $\R^{2n}$, the integral operator
$X \in S^\infty(\H)$ defined by
\begin{equation*}
(X\phi)(x) = \int_{\R^n} k(x,y) \phi(y) \, dy
\end{equation*}
is called a {\em Schwartz operator} on $\H$.  The class of Schwartz
operators on $\H$ is denoted by $\S(\H)$, and is a Fr\'echet space
with respect to an appropriate sequence of seminorms.  Moreover,
$\S(\H) \subseteq S^1(\H)$, and
$\alpha:\S(\H) \to \S(\R^{2n})$ and $W:\S(\R^{2n}) \to \S(\H)$ are
mutually inverse topological isomorphisms.  A continuous
linear functional on $\S(\H)$ is regarded as a distributional operator
on $\H$; and the aggregate of these is denoted $\S'(\H)$.
An element $X\in S^\infty(\H)$ defines a distributional operator $T_X$ on $\H$ by
$T_X(Y)=\tr(XY)$, $Y\in \S(\H)$.
Now, if $T$ is a tempered distribution on $\R^{2n}$, we define its Weyl
transform $W(T) \in S'(\H)$ by
\begin{equation*}
W(T)(Y)=T(\alpha(Y)^\sim),
\end{equation*}
where $\alpha(Y)^\sim(x,y)=\alpha(Y)(-x,-y)$.
The distributional Weyl transform is a topological isomorphism between
the space $\S'(\R^{2n})$ of tempered distributions on $\R^{2n}$ and the space
$\S'(\H)$.  Its inverse is
the distributional Fourier-Wigner transform, which is defined by
\begin{equation*}
\alpha(X)(\phi) = X(W(\phi^\sim)), \qquad X\in\S'(\H), \quad \phi\in\S(\R^{2n}).
\end{equation*}
\begin{thm}\label{T:LS}
Let $T$ be a compactly supported distribution on $\R^{2n}$.
Let $\hat{T}$ denote the Fourier transform of $T$. 
Then $W(T) \in S^p(\H)$ if and only if $\hat{T} \in L^p(\R^{2n})$
for $1 \leq p \leq \infty$.
Moreover, if $K$ is a compact set in $\R^{2n}$, then there exists a constant
$C_K$ such that
\begin{equation*}
C_K^{-1} \norm{\hat{T}}_p \leq
\norm{W(T)}_{S^p} \le C_K \norm{\hat{T}}_p,
\end{equation*}
whenever $\supp(T)\subseteq K$.
Furthermore, $W(T)$ is compact if and only if $\hat{T}$ vanishes at $\infty$.
\end{thm}

When $T$ is a Radon measure (i.e., when the order of $T$ is zero), this
theorem was proved by Luef and Samuelsen (see \cite[Theorem 1.4]{LS}).
Here we give a short proof of Theorem \ref{T:LS} based on
Tim Steger's observation that Fourier-Wigner transforms of
trace-class operators are locally the same as Fourier transforms
of integrable functions
(cf \cite{mkv-thesis}; see also \cite[Theorem 4.3 and Corollary 4.5]{rcr}).

In the proof of Theorem \ref{T:LS}, it is convenient to work with the
symplectic Fourier transform instead of the Fourier transform; this is
inconsequential because the former is obtained from the latter by
composing with a rotation of the variable, and hence the two are
comparable, in the sense that they belong to the Lebesgue spaces
simultaneously, and vanish at infinity simultaneously.
Recall that the {\em symplectic Fourier transform} of a function
$f\in L^1(\R^{2n})$ is the function on $\R^{2n}$ defined by
\begin{equation*}
\breve{f}(\xi, \eta)=
\int_{\R^n}\int_{\R^n}  f(x,y) e^{2\pi i(\xi \cdot y - \eta \cdot x)} \,dx\,dy.
\end{equation*}
More generally (see \cite{maillard}), if $T$ is a tempered distribution
on $\R^{2n}$, the symplectic Fourier transform of $T$ is given by
\begin{equation*}
\breve{T}(\phi)= T((\breve\phi)^\sim),
\quad \phi\in \mathcal{S}(\R^{2n}).
\end{equation*}

In \cite{rcr}, $S^1(\H)$ was regarded as an $L^1(\R^{2n})$-module under
the action
\begin{equation*}
f\cdot X = \int_{\R^n}\int_{\R^n} f(x,y) \p(x,y)X\p(x,y)^{-1} \,dx \,dy,
\qquad X\in S^1(\H), \quad f\in L^1(\R^{2n}),
\end{equation*}
and it was shown that $\alpha(f\cdot X)=\breve{f}\alpha(X)$.  We may
also regard $S'(\H)$ as a $\S(\R^{2n})$-module under the action
\begin{equation*}
(f\cdot X)(Y) = X(f^\sim \cdot Y),
\qquad X\in \S'(\H), \quad Y\in \S(\H), \quad f\in \S(\R^{2n}).
\end{equation*}
This coincides with the $L^1(\R^{2n})$-module structure
when both are defined,
and it is easy to check that $\alpha(f\cdot X)=\breve{f}\alpha(X)$
in this setting.

Let $g(x,y)= e^{-\frac{\pi}{2}(\abs{x}^2+\abs{y}^2)}$, $(x,y) \in \R^{2n}$.
For $X \in \S'(\H)$, let $\beta(X)=g\alpha(X)$ and
$\breve{\beta}(X)=\beta(X)\breve{}$.
According to \cite[Theorem 4.3]{rcr},
$\norm{\breve{\beta}(X)}_1 \leq \norm{X}_{S^1}$.
By Plancherel theorems for the Fourier and Weyl transforms
(see e.g., \cite{folland}), it follows that
\begin{equation*}
\norm{\breve{\beta}(X)}_2
=
\norm{\beta(X)}_2
\leq
\norm{g}_{\infty}\norm{\alpha(X)}_2
=
\norm{\alpha(X)}_2
=
\norm{X}_{S^2}.
\end{equation*}
By the Calderon-Lions interpolation theorem 
(\cite[Theorem 9.20, Example 1, Proposition 8]{RS}),
it follows that
\begin{equation}\label{E:beta-p}
\norm{\breve{\beta}(X)}_p \leq \norm{X}_{S^p},
\qquad X\in S^p(\H), \quad 1\le p \le 2.
\end{equation}

For $T \in \S'(\R^{2n})$, let $\Gamma(T) =W(g\breve{T})$.  Observe that
$\norm{W(g)}_{S^1} =1$, and hence by \cite[Lemma 3.10]{rcr},
$\norm{\Gamma(f)}_{S^1} \leq \norm{f}_1$ for $f\in L^1(\R^{2n})$.
An argument analogous to the one used to derive Equation (\ref{E:beta-p})
gives
\begin{equation}\label{E:gamma-p}
\norm{\Gamma(f)}_{S^p} \leq \norm{f}_p,
\qquad f\in L^p(\R^{2n}), \quad 1\le p\le 2.
\end{equation}

Now suppose $p>2$ and $X\in S^p(\H)$.  Let $p'$ be conjugate to $p$.  Then
$p'\in [1,2)$.  Regard $X$ as a distributional operator; then
for $\phi\in\S(\R^{2n})$, we have
\begin{equation*}
\abs{(\breve{\beta}(X))(\phi)}
= \abs{X(\Gamma(\phi))}
\le \norm{X}_{S^p} \norm{\Gamma(\phi)}_{S^{p'}}
\le \norm{X}_{S^p} \norm{\phi}_{p'}
\end{equation*}
by Equation (\ref{E:gamma-p}).
Therefore $\breve{\beta}(X) \in L^p(\R^{2n})$
and $\norm{\breve{\beta}}\le \norm{X}_{S^p}$.
Therefore Equation (\ref{E:beta-p}) holds for all $p\in [1,\infty]$.
Similarly, Equation (\ref{E:gamma-p}) also holds for all $p\in [1,\infty]$.



We are now ready to prove Theorem \ref{T:LS}.
Suppose $T$ is supported in $K$, and $\breve{T} \in L^p(\R^{2n})$.
Let $f$ be a compactly supported smooth function which is identically one
on $K$, and put $h(x,y)=e^{\frac{\pi}{2}(\abs{x}^2+\abs{y}^2)} f(x,y)$.  Then
$\breve{h}$ is integrable; put $C_K=\norm{\breve{h}}_1$.
Let $T_1= hT$.
By Young's inequality,
$\breve{T_1}=\breve{h}*\breve{T} \in L^p(\R^{2n})$, and
$\norm{\breve{T_1}}_p \leq C_K \norm{\breve{T}}_p$.
Therefore, by Equation \eqref{E:gamma-p}, $\Gamma(\breve{T_1})\in S^p(\H)$ and
$\norm{\Gamma(\breve{T_1})}_{S^p}\leq C_K \norm{\breve{T}}_p$.
However
\begin{equation*}
\Gamma(\breve{T_1})=W(gT_1)=W(ghT)=W(T).
\end{equation*}
Therefore $W(T) \in S^p(\H)$, and
$\norm{W(T)}_{S^p}\leq C_K \norm{\breve{T}}_p$.

Conversely, suppose $T$ is supported in $K$, and $W(T)\in S^p(\H)$.
Let $Z=\breve{h}\cdot W(T)$.
Since $\rho(x,y)$ is unitary for all $(x,y)\in\R^{2n}$, it follows that
\begin{equation*}
\norm{Z}_{S^p} \le
\int_{\R^n}\int_{\R^n} \abs{\breve{h}(x,y)} 
                        \norm{\rho(x,y)W(T)\rho(x,y)^{-1}}_{S^p} \,dx\,dy =
\norm{\breve{h}}_1 \norm{W(T)}_{S^p}.
\end{equation*}
Therefore $Z\in S^p(\H)$, and $\norm{Z}_{S^p} \leq C_K \norm{W(T)}_{S^p}$.
Therefore, by Equation \eqref{E:beta-p}, $\breve{\beta}(Z)\in L^p(\R^{2n})$ and
$\norm{\breve{\beta}(Z)}_p\leq C_K \norm{W(T)}_{S^p}$.
However,
\begin{equation*}
\breve{\beta}(Z)=
(g\alpha(Z))\breve{}=(ghT)\breve{}
=\breve{T}.
\end{equation*}
Therefore $\breve{T} \in L^p(\R^{2n})$, and
$$
C_K^{-1} \norm{\breve{T}}_p \le \norm{W(T)}_{S^p}.
$$

Now suppose $\breve{T} \in \mathcal{C}_0(\R^{2n})$.  
Let $B$ denote the closed unit ball centered at $0$ in $\R^{2n}$ and
put $K=\supp(T)+B$ (Minkowski sum).  Then $K$ is compact.  Let $\p$
be a non-negative smooth function supported in $B$ with $\int \p = 1$.
For each $r>1$, let $\p_r$ be defined by $\p_r(x,y)=r^{2n}\p(rx, ry)$.
Observe that for each $r>0$, $\breve\p_r\in L^1(\R^{2n})$,
$\norm{\breve\p_r}_\infty=1$, and $\lim_{r\to\infty} \breve\p_r = 1$ uniformly
on compact sets.  Let $\e>0$. Since $\breve{T}$ vanishes at $\infty$,
it follows that for sufficiently large $r$,
$\norm{\breve\p_r \breve{T} - \breve{T}}_\infty < \e/C_K$.  Put
$Y=W(\p_r * T)$.  Since $\p_r * T \in \S(\R^{2n})$, it follows that
$Y\in \S(\H)$, and hence compact.
Since $\p_r * T - T$ is supported in $K$, it follows that
\begin{equation*}
\norm{Y-W(T)}_{S^\infty}
  = \norm{W(\p_r * T -T)}_{S^\infty}
\le C_K \norm{\breve\p_r \breve{T} - \breve{T}}_\infty
  < \e.
\end{equation*}
Since the set of compact operators is closed in $S^\infty(\H)$, it follows
that $W(T)$ is compact.

Conversely, assume that
$W(T)$ is compact.
Observe that for each $r>0$, $W(\rho_r) \in S^1(\H)$, 
$\norm{W(\p_r)}_{S^\infty} \le 1$, and
$W(\rho_r)$ converges in the strong operator topology to the identity operator
as $r\to \infty$.
Let $\e>0$. 
Since $W(T)$ may be approximated by a finite rank operator,
it follows that for sufficiently large $r$,
$\norm{W(\rho_r)W(T)-W(T)}_{S^\infty} <\e/C_K$.
Let $f=\rho_r \twisted T$ be the twisted convolution of $\p_r$ with $T$
(see \cite{maillard}).
Then $\breve{f} \in \S(\R^{2n}) \subseteq \mathcal{C}_0(\R^{2n})$.
Since $\p_r \twisted T - T$ is supported in $K$, it follows that
\begin{equation*}
\norm{\breve{f}-\breve{T}}_{\infty} \le
C_K \norm{W(\p_r \twisted T -T)}_{S^\infty} =
C_K \norm{W(\p_r)W(T)-W(T)}_{S^\infty} <
\e.
\end{equation*}
Since $\mathcal{C}_0(\R^{2n})$ is closed in $L^\infty(\R^{2n})$,
it follows that $\breve{T} \in \mathcal{C}_0(\R^{2n})$.

\bibliographystyle{amsplain}
\bibliography{v8-wtcsd}

\end{document}